\documentclass[sn-mathphys-num]{sn-jnl}% Math and Physical Sciences Numbered Reference Style 
\usepackage{graphicx}
\usepackage{float}
\usepackage{amsmath}
\usepackage{amsfonts}
\graphicspath{{Figures/}} 
\usepackage{subcaption}
\begin{document}
\title{Delay Infectivity and Delay Recovery SIR model}
\author{C. N. Angstmann, S-J. M. Burney, A. V. McGann, Z. Xu}
\maketitle

\begin{abstract}
	
The governing equations for an SIR model with delay terms in both the infectivity and recovery rate of the disease are derived. 
The equations are derived by modelling the dynamics as a continuous time random walk, where individuals move between the classic SIR compartments. Delays are introduced through an appropriate choice of distribution for the infectivity and recovery processes. This approach ensures the physicality of the model and provides insight into the underlying dynamics of existing SIR models with infectivity delays. We show that a delay in the infectivity term arises from taking the infectivity to be a hazard function of the disease recovery rate, and a delay in the recovery arises from taking the corresponding survival function as a delay exponential distribution. The delays in the derived model can represent a latency period or incubation effect. Compared to the standard SIR model, this framework allows for a broader range of dynamical behaviors while remaining highly tractable.
\end{abstract}

\section{Introduction}
The spread of infectious diseases through populations have been modelled by SIR models, since they were first introduced by Kermack and McKendrick in 1927 \cite{KM1927}. In this model, the population is split into three compartments - those susceptible (S), infective (I) and recovered (R) from the infection. Individuals in the population move through the compartments. The time-evolution of the population of the compartments is represented by a set of three coupled ordinary differential equations (ODEs). In the intervening years many extensions have been proposed to this model to account for further dynamics, including incorporating `age-of-infection' effects \cite{KM1932}, stratifying the population based on sex \cite{H2000,M2003} or the inclusion of additional compartments, such as the SEIR model which includes an `Exposed' subclass of the population \cite{HT1980}. 

An increasingly popular way to incorporate `time-since-infection' effects is through the generalisation of the ODEs to a system of delay differential equations (DDEs) \cite{RHRC2018}. These models can better reflect particular disease dynamics, when accounting for effects such as an incubation time or a latency period \cite{C1979,CD1996}. There have been many generalised SIR models proposed that incorporate a time-delay on the infectivity rate, within the `force of infection' term \cite{TNL2007,MST2004,BT1995,HD2000,Liu2015}. Delays have also been considered in more general epidemiological models, such as SEIR and SEIRS models \cite{CD1996, YL2006,HTMW2010}. A delay on the re-susceptibility of individuals after a temporary immunity from the infection has also been considered \cite{CD1996}. The delays in the infectivity rate have been incorporated in multiple ways. In some model formulations the delay only affects the infective population in the interaction term \cite{TNL2007,MST2004,BT1995,YL2006,HTMW2010,Liu2015}. However, in Cooke and Van den Driessche \cite{CD1996}, an SEIRS model with two delays is established. In this model, both the susceptible and infective populations are subject to a delay in the bilinear interaction term, representing an infection incubation period. The delay requires the use of the `Exposed' class, and the SEIRS model reduces to an SIRS model if the incubation delay is taken to be zero. There has also been some inconsistency in the interpretation of model parameters across the literature.

In order to gain further insight on the processes governing an SIR model with delays, the model can be derived from the underlying stochastic process. This stochastic process describes the evolution of the infection through the population and ensures the resulting DDE is consistent with the governing equation of the process. The resulting SIR model features both a delay on the infectivity and recovery rates.
We begin by deriving an SIR model from a continuous time random walk (CTRW) \cite{MW1965}, and consider an arbitrary infectivity and recovery rate \cite{AHM2017}. This approach has previously been used to show the conditions under which fractional derivatives arise in generalised SIR models \cite{AHM2016,AHM2016fi,AEHMMN2017,AEHMMN2021,AHM2017}. 
We consider the necessary forms of the infectivity and recovery in the underlying stochastic process that lead to time-delay terms in the SIR model. For the recovery, we require the time spent in the infected compartment to follow a delay exponential distribution \cite{ABHJX2023,ABHHMX2024}. This leads to a delay in the recovery term. We also consider a form of the infectivity or `force-of-infection' which leads to delay in the infectivity term. Both of these effects can be taken together to produce a delay infectivity and delay recovery SIR model.
By deriving the governing equations from a stochastic process, the resulting model is comprised of intuitive model parameters.
Additionally, this method allows us to establish critical values on the resulting delay terms. Hence, we are able to provide further insight on existing SIR models with delay infectivity terms.

In Section \ref{sec:deriv} we derive a model with delays on the infectivity and recovery and consider the critical values on the time-delays and steady states of the system. In Section \ref{sec_specialcases} we consider reductions to the delay infectivity SIR and delay recovery SIR models. Examples are shown in Section \ref{sec:examples} and we conclude with a discussion in Section \ref{sec:summary}.

\section{Derivation}\label{sec:deriv}
In order to incorporate a delay infectivity and a delay recovery into an SIR model with vital dynamics, we begin by taking the general set of master equations from a CTRW. In the CTRW, we consider individuals entering a compartment, waiting there for a random amount of time before leaving to another compartment until they leave the system through death. The length of time an individual spends in a compartment is drawn from a waiting time distribution. We consider an arbitrary infectivity and recovery rate as in Angstman et al. \cite{AHM2017} where integro-differential equations govern the populations of the Susceptible (S), Infective (I) and Recovery (R) compartments. A full derivation of the model is provided in Appendix \ref{sec_initial_cond}. 

The master equations for the time evolution of the epidemic are:
\begin{align}
\frac{dS(t)}{dt}&=\lambda(t)-\omega(t) S(t) \theta(t,0)\int_{0}^{t}K_I(t-t')\frac{I(t')}{\theta(t',0)}dt'-\gamma(t)S(t),\label{eq_genS_i}\\
\frac{dI(t)}{dt}&=\omega(t) S(t) \theta(t,0)\int_{0}^{t}K_I(t-t')\frac{I(t')}{\theta(t',0)}dt'
-\theta(t,0)\int_{0}^{t}K_R(t-t')\frac{I(t')}{\theta(t',0)}dt'-\gamma(t)I(t),\label{eq_genI_i}\\
\frac{dR(t)}{dt}&=\theta(t,0)\int_{0}^{t}K_R(t-t')\frac{I(t')}{\theta(t',0)}dt'-\gamma(t)R(t)\label{eq_genR_i}.
\end{align}
In this set of equations, $\lambda(t)>0$ is the birth rate and $\gamma(t)>0$ is the death rate per capita. The environmental infectivity rate is $\omega(t)$ and the probability of surviving the death process from time $t'$ to $t$ is captured by $\theta(t,t')$. The initial conditions of the infectivity compartment are taken such that $I(t)=0$ for $t<0$.

Individuals may only enter the Infective compartment from the Susceptible compartment, hence there is a corresponding decrease in the number of individuals in the Susceptible compartment. Similarly, the individuals who leave the Infective compartment through recovery, correspond with the flux into the Recovery compartment. 

The infectivity ($K_I$) and recovery ($K_R$) memory kernels are the result of taking Laplace transforms to enable us to write governing equations.
The memory kernel of the recovery function is
\begin{equation}\label{eq:K_recovery}
K_R(t)=\mathcal{L}_s^{-1}\left\{\frac{\mathcal{L}_t\left\{\psi(t)\right\}}{\mathcal{L}_t\left\{\phi(t)\right\}}\right\}.
\end{equation}
Here $\phi(t)$ is the probability of not recovering from the infected state after time $t$, and $\psi(t)$ is defined as the corresponding waiting time probability density function, hence,
\begin{equation}
\psi(t)=-\frac{d\phi(t)}{dt}.
\end{equation}
The memory kernel of the infectivity is defined as
\begin{equation}\label{eq:K_infectivity}
K_I(t)=\mathcal{L}_s^{-1}\left\{\frac{\mathcal{L}_t\left\{\rho(t)\phi(t)\right\}}{\mathcal{L}_t\left\{\phi(t)\right\}}\right\},
\end{equation}
where $\rho(t)$ is the age-of-infection dependent infectivity rate. If $\phi(t)$ is an exponential and $\rho(t)$ is a constant, the standard SIR model is recovered. 

To incorporate a delay into the infectivity and recovery terms, we will choose waiting time distribution, $\phi(t)$, and infectivity rate, $\rho(t)$, such that the convolution integrals will induce delays.
To obtain a delayed recovery term, we take the recovery survival function, $\phi(t)$, to be the survival function of the delay exponential distribution \cite{ABHHMX2024}. 
The delay exponential function is defined by the power series \cite{ABHJX2023},
\begin{equation} \label{def_dexp}
	\mathrm{dexp}(-\mu t;-\mu \tau_2)=\sum_{n=0}^\infty  \frac{(-\mu)^n(t-n\tau_2)^n}{\Gamma(n+1)}\Theta\left(\frac{t}{\tau_2}-n\right),\quad \frac{t}{\tau_2}\in  \mathbb{R},
\end{equation}
and the Heaviside function is defined by,
\begin{equation}
	\Theta(y)=\begin{cases} 0 & y<0,\\ 1& y\ge 0. \end{cases}
\end{equation}
Eq. \eqref{def_dexp} is only a valid survival function when $0 \leq \mu\tau_2 \leq e^{-1}$. Hence, we can express $\phi(t)$ as,
\begin{equation}\label{eq:phi}
\phi(t)=\text{dexp}(-\mu t;-\mu \tau_2),
\end{equation}
with $0 \leq \mu\tau_2 \leq e^{-1}$.
Here, $\tau_2$ represents a constant delay and $\mu^{-1}$ is the mean of the delay exponential distribution. 

\begin{figure}[H]
	\begin{center}
		\includegraphics[width=0.85\textwidth]{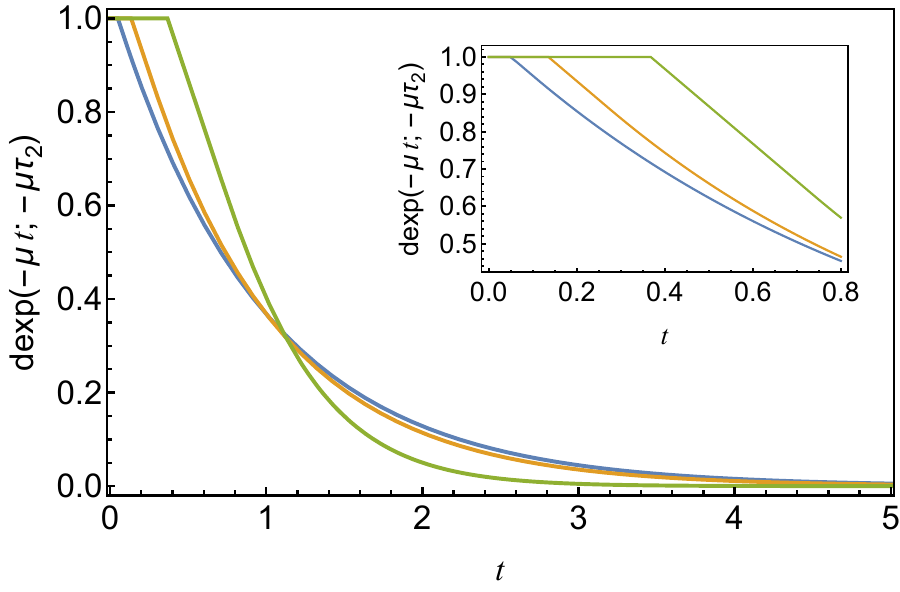}
		\caption{Plot of the delay exponential function, Eq. (\ref{def_dexp}). The green, orange and blue curves correspond to $\tau_2=e^{-1}, e^{-2}$ and $e^{-3}$, respectively, with $\mu=1$.}
		\label{fig:dexp}
	\end{center}
\end{figure}

The dynamics of Eq. (\ref{def_dexp}) are illustrated in Figure \ref{fig:dexp} for three different delay values such that $\mu\tau_2\in [0,e^{-1}]$. 

The Laplace transform of Eq. \eqref{eq:phi} is
\begin{equation}
\mathcal{L}_t\left\{\phi(t)\right\}=\frac{1}{s+\mu e^{-s\tau}},
\end{equation}
with the memory kernel,
\begin{equation}\label{eq:K_R}
K_R(t)=\mu \delta(t-\tau_2),
\end{equation}
where $\delta$ is the Dirac delta generalised function. 

The recovery convolution can then be explicitly obtained via the sifting property of the Dirac delta,
\begin{align}
\int_{0}^{t}K_R(t-t')\frac{I(t')}{\theta(t',0)}dt'&=\mu \int_{0}^{t}\delta(t-t'-\tau_2)\frac{I(t')}{\theta(t',0)}dt'\\
&=\mu \frac{I(t-\tau_2)}{\theta(t-\tau_2,0)}\Theta(t-\tau_2).\label{eq:K_rec}
\end{align}
Given initial conditions $I(t)=0$ for $t<0$, the equation can be written without the Heaviside function.

\begin{figure}[h!]
	\begin{center}
	\includegraphics[width=0.8\textwidth]{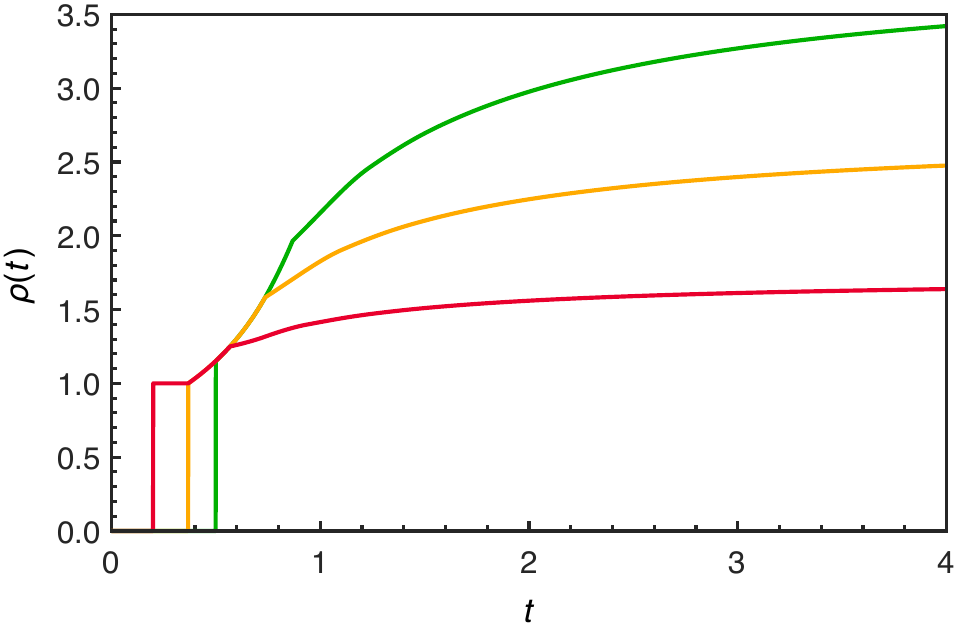}
	\caption{Plot of the infectivity rate, $\rho(t)$, with $\mu=1$ and $\tau_2=e^{-1}$. The red, orange and green curves correspond to $\tau_1= 0.2, e^{-1}$ and $0.5$, respectively.}
	\label{fig:dexp_rho}
\end{center}
\end{figure}

Now, we consider the conditions for a delay to be present in the infectivity term.  To force a delay term into the infectivity, we take $\rho(t)$ to be,
\begin{equation}\label{eq:rho}
\rho(t)=\frac{\Theta(t-\tau_1)\phi(t-\tau_1)}{\phi(t)}.
\end{equation}
This choice of $\rho(t)$ will lead to a time-delay in the infectivity, regardless of the $\phi(t)$ taken. When $\tau_1=\tau_2$, $\rho(t)$ is a hazard function of the chosen $\phi(t)$ distribution. The dynamics of the infectivity of $\rho(t)$ with a delay exponential recovery survival $\phi(t)$ is shown in Figure \ref{fig:dexp_rho}. We consider three cases of $\tau_1$ in this figure. In each case the infectivity begins at zero. The larger the value of $\tau_1$ the longer the infectivity stays at zero before `switching on'. Larger values of $\tau_1$ result in an increased infectivity once the infectivity `switches on'. Note that $\rho(t)$ has a $\tau_2$ dependence within it due to its $\phi(t)$ dependence.

The infectivity kernel can then be written as,
\begin{align}
K_I(t)&=\mathcal{L}_s^{-1}\left\{\frac{\mathcal{L}_t\left\{\Theta(t-\tau_1)\phi(t-\tau_1)\right\}}{\mathcal{L}_t\left\{\phi(t)\right\}}\right\}.
\end{align}
Similar to the recovery kernel, this simplifies to,
\begin{equation}
K_I(t)=\delta(t-\tau_1).
\end{equation}
The infectivity convolution can then be obtained similarly to Eq. \eqref{eq:K_R}, hence
\begin{align}
\int_{0}^{t}K_I(t-t')\frac{I(t')}{\theta(t',0)}dt'&=\mathcal{L}_s^{-1}\left\{\mathcal{L}_t\left\{K_I(t)\right\}\mathcal{L}_t\left\{\frac{I(t)}{\theta(t,0)}\right\}\right\}\\
&=\mathcal{L}_s^{-1}\left\{ e^{-s\tau_1}\mathcal{L}_t\left\{\frac{I(t)}{\theta(t,0)}\right\}\right\}\\
&=\frac{I(t-\tau_1)}{\theta(t-\tau_1,0)}\Theta(t-\tau_1).\label{eq:K_inf}
\end{align}
As with the recovery delay equation, Eq. \eqref{eq:K_rec}, the Heaviside function can be dropped due to the initial conditions. Taking these initial conditions and substituting Eqs. \eqref{eq:K_rec} and \eqref{eq:K_inf} into Eqs. \eqref{eq_genS_i}, \eqref{eq_genI_i} and \eqref{eq_genR_i}, we obtain the governing equations:
\begin{align}
\frac{dS(t)}{
dt}=&\lambda(t)-\omega(t) S(t) \theta(t,0)\frac{I(t-\tau_1)}{\theta(t-\tau_1,0)}-\gamma(t)S(t),\label{eq_genS_delay}\\
\frac{dI(t)}{dt}=&\omega(t) S(t) \theta(t,0)\frac{I(t-\tau_1)}{\theta(t-\tau_1,0)} -\mu\theta(t,0)\frac{I(t-\tau_2)}{\theta(t-\tau_2,0)}-\gamma(t)I(t),\label{eq_genI_delay}\\
\frac{dR(t)}{dt}=&\mu\theta(t,0)\frac{I(t-\tau_2)}{\theta(t-\tau_2,0)}-\gamma(t)R(t)\label{eq_genR_delay}.
\end{align}

\begin{figure}[H]
	\begin{center}
		\includegraphics[width=0.9\linewidth]{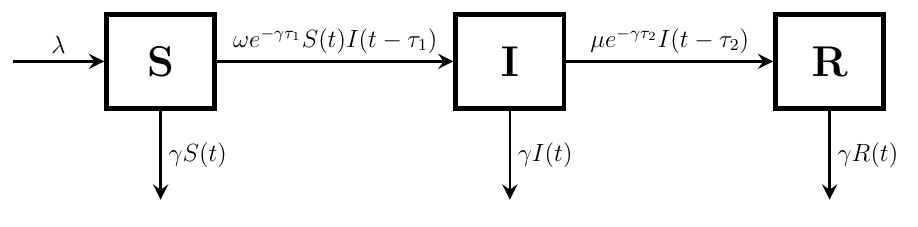}
		\caption{An SIR model with constant vital dynamics and a time-delay on the infectivity and recovery terms.}
		\vspace{-1.5\baselineskip}
		\label{fig:comp_model}
	\end{center}
\end{figure}

By taking the birth, death and infectivity rates to be constant, $\lambda(t)=\lambda$, $\gamma(t)=\gamma$ and $\omega(t)=\omega$ respectively, we can simplify the equations to:
\begin{align}
\frac{dS(t)}{dt}=&\lambda-\omega e^{-\gamma\tau_1} S(t)  I(t-\tau_1)-\gamma S(t),\label{eq_genS_delay_simp}\\
\frac{dI(t)}{dt}=&\omega e^{-\gamma\tau_1}S(t) I(t-\tau_1)-\mu e^{-\gamma\tau_2} I(t-\tau_2)-\gamma I(t),\label{eq_genI_delay_simp}\\
\frac{dR(t)}{dt}=&\mu e^{-\gamma\tau_2} I(t-\tau_2)-\gamma R(t)\label{eq_genR_delay_simp}.
\end{align}
We will consider this set of simplified equations for the remainder of the paper. A representation of the movement through compartments is shown in Figure \ref{fig:comp_model}. 

In this set of equations the parameters of the model retain their standard dimensions, as well as their interpretation. In some of the previous SIR models with time-delays, there is no accounting for the $e^{-\gamma\tau_1}$ survival term in the infectivity term. While this does not change the correctness of these previous models, it does require a different interpretation for the remaining infectivity term, $\omega$. It also shows that if a model considers a different size delay, the $\omega e^{-\gamma\tau_1}$ term will be rescaled. This reasoning holds for the $\tau_2$ survival term as well.

\subsection{Critical Values of the Time-Delays}
We now turn our attention to the critical values of delay parameters, $\tau_1$, $\tau_2$ and $\mu$. In Section \ref{sec:deriv}, we discussed that the delay exponential function is only a probability distribution when $0 \leq \mu\tau_2 \leq e^{-1}$. Hence, if we consider a large $\tau_2$ then $\mu$ must be adequately small. However, there is no similar restriction on $\tau_1$, as the infectivity is not governed by a probability function. The infectivity, $\rho(t)$, is only required to be non-negative. As $\rho(t)$ is defined by Eq. \eqref{eq:rho}, it will be non-negative for all $\tau_1\geq 0$, as it is composed of non-negative functions.

\subsection{Steady States}\label{sec_steadystates}
It is straightforward to find the steady states of the model defined by Eqs. \eqref{eq_genS_delay_simp}, \eqref{eq_genI_delay_simp} and \eqref{eq_genR_delay_simp}. We will define the steady state to be $(S^*, I^*, R^*)$ where,
\begin{equation}
\lim_{t\rightarrow \infty} S(t)=S^*,\hspace{30pt}\lim_{t\rightarrow \infty} I(t)=I^*,\hspace{30pt}\lim_{t\rightarrow \infty} R(t)=R^*.
\end{equation}
These steady state values will satisfy the equations:
\begin{align}
0=&\lambda-e^{-\gamma\tau_1} \omega S^* I^*-\gamma S^*,\label{eq_genS_steady}\\
0=&e^{-\gamma\tau_1} \omega S^*I^*-e^{-\gamma\tau_2}\mu  I^*-\gamma I^*,\label{eq_genI_steady}\\
0=&e^{-\gamma\tau_2} \mu I^*-\gamma R^*\label{eq_genR_steady}.
\end{align}
The disease-free steady state is:
\begin{equation}\label{eq:ss_free}
S^*=\frac{\lambda}{\gamma},\hspace{30pt}I^*=0,\hspace{30pt}R^*=0.
\end{equation}
The endemic steady state is:
\begin{equation} \label{eq:ss_endemic}
S^*=\frac{\mu e^{-\gamma \tau_2}+\gamma}{\omega e^{-\gamma \tau_1}},\hspace{12pt}I^*=\frac{\lambda}{\mu e^{-\gamma \tau_2}+\gamma}-\frac{\gamma}{\omega e^{-\gamma \tau_1}},\hspace{12pt}R^*=\frac{\lambda \mu e^{-\gamma \tau_2}}{\gamma\mu e^{-\gamma \tau_2}+\gamma^2}-\frac{\mu e^{-\gamma \tau_2}}{\omega e^{-\gamma \tau_1}}.
\end{equation}
Note that the endemic steady state exists only if
\begin{equation}\label{eq_endemic_cond}
\lambda\omega e^{-\gamma \tau_1}>\gamma\left(\mu e^{-\gamma \tau_2}+\gamma\right).
\end{equation}

\section{Reductions} \label{sec_specialcases}
The general delay infectivity and delay recovery SIR equations can be reduced to simpler SIR models. An SIR model with only a delay on the infectivity or recovery is produced by setting one of the delays to zero. 
This enables us to compare existing SIR models with time-delays on the infectivity to our reduced model, a delay infectivity SIR model, when $\tau_2=0$ and $\tau_1>0$. Note that by setting $\tau_1=\tau_2=0$ in the delay infectivity and delay recovery SIR model, the standard SIR model is recovered.
This shows that our delay infectivity and recovery SIR model is consistent with the standard SIR model. In this section, we present the reduced models.

\subsection{Delay infectivity SIR}\label{sec_disir}
To obtain a delay infectivity SIR model, we set $\tau_2=0$. Hence, the set of governing equations for the SIR model with a delay infectivity term are:
\begin{align}
\frac{dS(t)}{dt}=&\lambda-\omega S(t) e^{-\gamma\tau_1} I(t-\tau_1)-\gamma S(t),\label{eq_genS_delayinf}\\
\frac{dI(t)}{dt}=&\omega S(t)e^{-\gamma\tau_1}I(t-\tau_1)-\mu I(t)-\gamma I(t),\label{eq_genI_delayinf}\\
\frac{dR(t)}{dt}=&\mu I(t)-\gamma R(t)\label{eq_genR_delayinf}.
\end{align}
When $\tau_2=0$ then Eq. \eqref{eq:phi} shows $\phi(t)$ to be exponentially distributed, hence $\phi(t)=e^{-\mu t}$. As the infectivity, Eq. \eqref{eq:rho}, is dependent on the recovery waiting time, we can identify the $\rho(t)$ that leads to the existing SIR models with infectivity delays, 
\begin{equation}\label{eq:rho_disir}
\rho(t)=\Theta(t-\tau_1)e^{\mu\tau_1}.
\end{equation}
Hence, SIR models with a delay infectivity rate are underpinned by no infectivity until an individual has been infected for a $\tau_1$ length of time, and then a constant rate. A substitution of $\tau_2=0$ into the steady states in Eqs. \eqref{eq:ss_free} and \eqref{eq:ss_endemic} gives us the steady states for this model.
\subsection{Delay recovery SIR}\label{sec_drsir}
To recover the delay recovery SIR model, we set $\tau_1=0$. Hence the set of governing equations for this model are:
\begin{align}
\frac{dS(t)}{dt}=&\lambda-\omega S(t) I(t)-\gamma S(t),\label{eq_genS_delayrec}\\
\frac{dI(t)}{dt}=&\omega S(t)I(t)-\mu e^{-\gamma\tau_2} I(t-\tau_2)-\gamma I(t),\label{eq_genI_delayrec}\\
\frac{dR(t)}{dt}=&\mu e^{-\gamma\tau_2} I(t-\tau_2)-\gamma R(t)\label{eq_genR_delayrec}.
\end{align}
When $\tau_1=0$, the infectivity, Eq. \eqref{eq:rho}, becomes $\rho(t)=1$ for $t\geq0$. Meanwhile, the infection recovery waiting time remains a delay exponential distribution. A substitution of $\tau_1=0$ into the steady states in Eqs. \eqref{eq:ss_free} and \eqref{eq:ss_endemic} gives us the steady states for this model.

\section{Results} \label{sec:examples}
In this section, we explore the effects of different delay parameters on the delay infectivity and delay recovery SIR model. We find that changes to the delay parameters lead to significant impacts in the short term dynamics of the model but cause minimal impacts in the long term dynamics.
We will consider the effect of varying the delay parameters, $\tau_1$, $\tau_2$ as well as the timescale parameter, $\mu$. We note that the variation in $\tau_1$ leads to changes in the infective rate, Eq. \eqref{eq:rho}. While variations in $\tau_2$ and $\mu$ affects both the infectivity rate, Eq. \eqref{eq:rho}, and recovery survival function, Eq. \eqref{eq:phi}.
For this study, we have taken, $\lambda=0.5$, $\gamma=0.001$ and $\omega=0.02$.

\begin{figure}[h!]
	\begin{center}
		\includegraphics[width=0.9\textwidth]{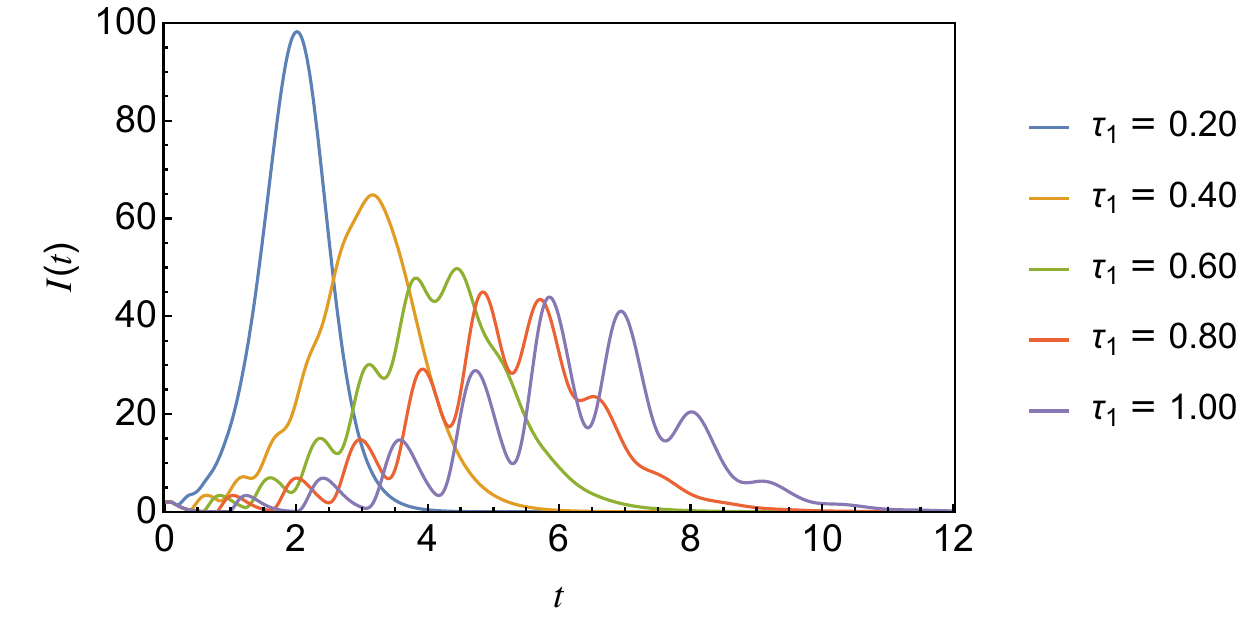}
		\caption{Plot of the Infective compartment when $\tau_1$ ranges from $0.5$ and $2.5$, with $\tau_2=0.1$ and $\mu=10e^{-1}$. The birth rate is $\lambda=0.5$, the death rate is $\gamma=0.001$, the infectivity rate is $\omega=0.02$, with initial values $S(0)=498$, $I(0)=2$ and $R(0)=0$.}
		\label{fig:tau_1_var}
	\end{center}
\end{figure}

We begin by varying the infectivity delay, $\tau_1$. We set the recovery delay, $\tau_2=0.1$ and $\mu$ has been defined as,
\begin{equation} \label{eq:mu_def}
\mu=\frac{e^{-1}}{\tau_2},
\end{equation}
to ensure $0\leq\mu\tau_2\leq e^{-1}$.
The larger the infective delay, the longer the infectivity is zero before the infectivity `turns on'. The $\tau_1$ delay reduces the peak of the infective compartment, as compared to a standard SIR model with matching vital dynamics and constant infective rate. The impact of varying $\tau_1$ on the population of the Infective compartment can be seen in Figure \ref{fig:tau_1_var}. 

\begin{figure}[h!] 
	\begin{center}
		\includegraphics[width=0.9\textwidth]{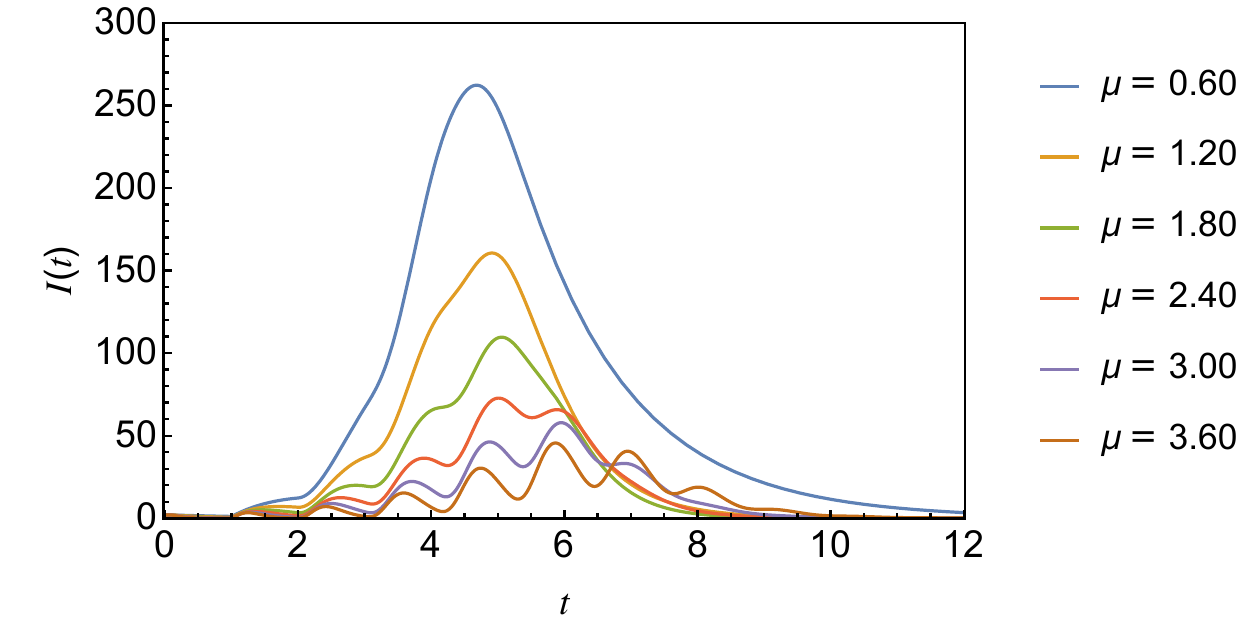}
		\caption{Plot of the Infective compartment with $\tau_1=1$ and $\tau_2=0.1$, with values for $\mu$ between $0.6$ and $3.6$. The birth rate is $\lambda=0.5$, the death rate is $\gamma=0.001$, the infectivity rate is $\omega=0.02$, with initial values $S(0)=498$, $I(0)=2$ and $R(0)=0$.}
		\label{fig:large_tau1_var_tau2}
	\end{center}
\end{figure}

\begin{figure}[H] 
	\begin{center}
		\includegraphics[width=0.9\textwidth]{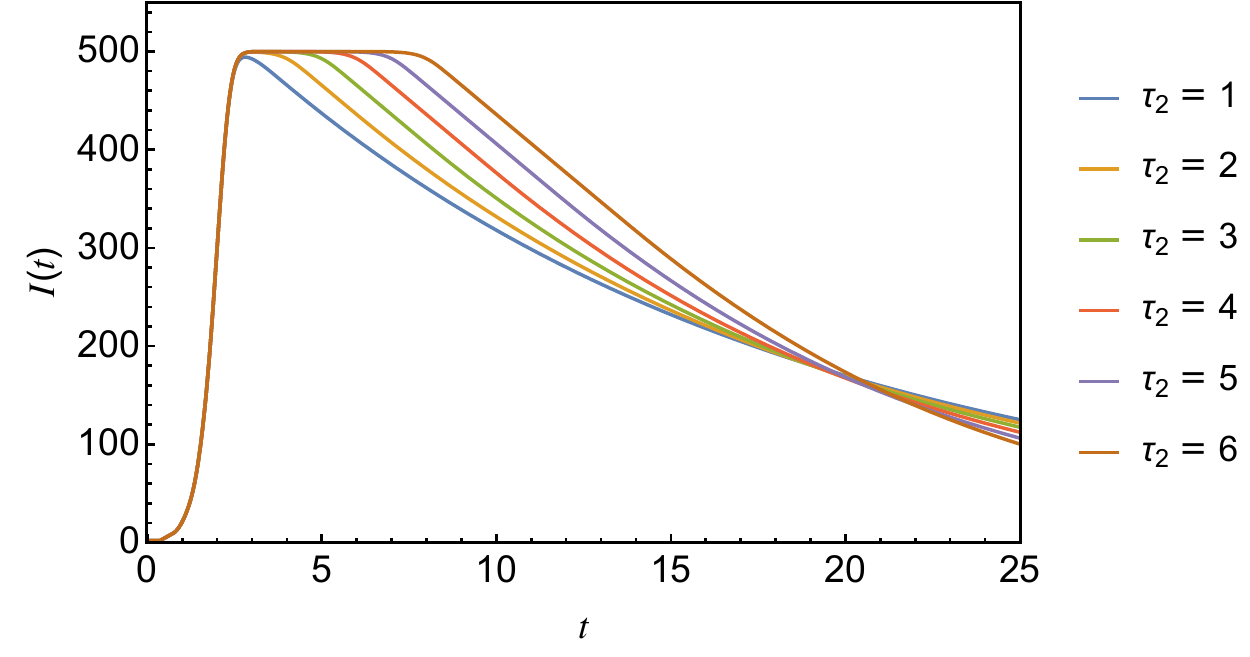}
		\caption{Plot of the Infective compartment when $\tau_1=e^{-1}$, $\mu=0.06$ and $\tau_2$ is between $1$ and $6$. The birth rate is $\lambda=0.5$, the death rate is $\gamma=0.001$, the infectivity rate is $\omega=0.02$, with initial values $S(0)=498$, $I(0)=2$ and $R(0)=0$.}
		\label{fig:tau_2_var_steady}
	\end{center}
\end{figure}

The larger the infective delay, the smaller the infective peak and the later the peak occurs. Of more note, the infective delay induces oscillations in the infective compartment, not observed in the standard SIR model. The oscillations occur at a period of the order of $\tau_1$. The oscillations persist longer for greater values of $\tau_1$. 
\newpage
Next, we consider the impact of varying $\mu$ on the Infective compartment. As $\mu$ is decreased, the oscillation effect is dampened. This result is presented in Figure \ref{fig:large_tau1_var_tau2}. In this figure, $\tau_1=1$ and $\tau_2=0.1$. Note that $\mu\tau_2\leq e^{-1}$ for all of the plot lines to ensure the Infective compartment remains non-negative. We see that as $\mu$ is decreased, a more standard SIR Infective compartment with a unimodal distribution is recovered, with an increased peak infectivity.
Next, we consider the effect of varying $\tau_2$. The oscillations persist when $\tau_2$ is increased, until $\tau_2$ is greater than $\tau_1$.
In Figure \ref{fig:tau_2_var_steady}, we set $\tau_1=e^{-1}$, $\mu=0.06$ and consider larger values for $\tau_2$. This leads to a sustained peak infection in the compartment. The larger the delay, $\tau_2$, the longer the peak infection number is sustained.

The sustained peak infectivity is an intuitive result given the vital dynamics and choice of  $\tau_2$. Under these conditions, the infective population grows until almost the entire population is infected. Infected individuals then have two ways to leave the infected compartment, either through recovery or death. The value $\tau_2$ traps individuals in the infected compartment for $\tau_2$ days before individuals have the ability to recover and with a much smaller death rate, the probability of dying is minimal. If the vital dynamics are increased, the steady peak infectivity occurs for less time. A different choice of vital dynamics would be considered in disease processes occurring over different magnitudes of time, i.e. the disease process of influenza where a typical infection lasts 6 to 8 days \cite{CK1983,PHS2016} would require a different timescale to a disease such as  human papillomavirus (HPV), where infections can take over a year to clear \cite{IDELB2007}.

Lastly, the reductions of the model from Subsections \ref{sec_disir} and \ref{sec_drsir} are considered. When $\tau_2=0$, we recover the delay infectivity SIR model from Subsection \ref{sec_disir}. In this model, $\mu$ is no longer bound by Eq. \eqref{eq:mu_def} and can be made arbitrarily large. Increasing either $\mu$ or $\tau_1$ introduces stronger oscillations into the system. When we consider the delay recovery SIR reduction, as in Subsection \ref{sec_drsir} with $\tau_1=0$, no oscillations occur. For larger values of $\tau_2$ we still observe a sustained infective peak, however larger choices of $\mu$ will dampen this impact.

Overall, the delay parameters produce substantial impacts to the short term dynamics of the epidemiological model. We have observed that for $\tau_1>\tau_2$ oscillations are produced, and for $\tau_2>\tau_1$ a steady peak is produced. However, when $\tau_1$ and $\tau_2$ are similar values their effects cancel out. The effects of both delays are also minimised, when the vital dynamics are increased.

An advantage of the delay infectivity and delay recovery SIR model is its computational tractability. The inclusion of the delay term in generalised SIR models is often done with the aim of capturing a `time-since-infection' or other historical effect. There has also been rising interest in capturing the impact of changed behaviour in the population in response to a known disease outbreak \cite{MD2013}. In some of these instances, the disease process may be more accurately represented by a set of integro-differential equations. However, there may be circumstances, where the accuracy can be sacrificed in order to be able to solve the system of equations. In such cases, the delay infectivity and delay recovery SIR model provides a way to incorporate past states as well as be able to solve the equations efficiently.

\section{Discussion}\label{sec:summary}
We have derived an SIR model with delays on both the infectivity and recovery terms by defining the SIR model from an underlying stochastic process. The delay terms are a consequence of our chosen infectivity and recovery functions. Our approach leads to better understanding of the mechanics at play when time-delays are added into SIR model equations.

We found that setting the recovery survival function to be a delay exponential distribution returned a delay in the recovery term. In order to incorporate a delay into the infectivity term, we defined the infectivity as a hazard function of delay exponential distributions. In the literature, there have been multiple SIR models with delays on the infectivity terms presented. By deriving our model from an underlying stochastic process, we can better understand the differences between the existing models. Our approach also provides insight into the parameters in the model, in particular showing that there is a survival term that is introduced when a delay is incorporated into the system. 

We established critical values of the delay on the recovery term, by examining the constraints on the delay exponential distribution. It was also observed that the recovery delay impacts the infectivity function as the infectivity is defined through the recovery density function. The only restriction on the delay on the infectivity term is that it must remain positive.

The steady states of the system were found as well as the condition on the endemic steady state existing. There remains the open question of the stability of these steady states, particularly under a range of initial conditions. This has been the focus of previous delay SIR and SEIR models \cite{C1979,MST2004,BT1995,TNL2007,HTMW2010,GOO2023}, however the existing models are dependent on only one delay term. There has been an approach put forth to consider the stability of first-order linear DDE systems with multiple delay terms \cite{G2010}.

The conditions under which the model reduced to having a single delay, either on the infectivity or the recovery, were established. We considered the implication of such a choice of function on the mechanisms governing the infection. The sets of equations to govern the delay infectivity SIR and delay recovery SIR models were provided.

Finally, we showed some of the ways the delay infectivity and delay recovery SIR model can differ from a classic SIR model. We examined the impact of varying all three parameters related to the delays. It was shown that when the delay in the infectivity is larger than the delay in the recovery, this induces oscillations into the Infective compartment. We showed that decreasing the timescale parameter, $\mu$, inhibits the magnitude of the oscillations. It was also found that when the delay in the recovery is larger than the delay in the infectivity, a steady infectivity peak is produced. 

In some instances, considering a set of integro-differential equations may serve to be a more accurate representation of the dynamics of a disease process. However, the results from our model show that through considering only two fixed time-delays, a simple change in parameter values can result in very different disease dynamics. Hence, there is potential that these dynamics may adequately approximate certain integro-differential SIR models. Additionally we have demonstrated that the delay infectivity and delay recovery SIR model is highly tractable. The tractability may, in some instances, provide enough benefit to choose this model over a more accurate set of integro-differential SIR equations.

\appendix
\section{Appendix} \label{sec_initial_cond}

We begin by defining the flux into the Infected compartment at time $t$. This will be defined as $q^+(I,t)$. The flux is constructed recursively as
\begin{equation}
q^+(I,t)=\int_{-\infty}^t \rho(t-t')\omega(t)S(t)\theta(t,t')\phi(t-t')q^+(I,t')dt'.
\end{equation}

Here, $\rho(t)$ is the age-of-infection dependent infectivity rate and $\omega(t)$ is the environmental dependent infectivity. Note that $\theta(t,t')$ is the probability of surviving the death process from $t'$ to $t$, and we assume it is of the form
\begin{equation}\label{eq:theta}
\theta(t,t')=\exp\left(-\int_{t'}^t\gamma(s)ds\right).
\end{equation}
Hence, this function obeys the semi-group property.
The probability of surviving the transition into the recovery compartment from $t'$ to $t$ is defined as $\phi(t-t')$. The recovery rate, $\psi(t)$, is subsequently defined as
\begin{equation}
\psi(t)=-\frac{d\phi(t)}{dt}.
\end{equation} 
We also define the flux, for $t<0$, as
\begin{equation}\label{eq:pre_flux}
q^+(I,t)=\frac{i(-t,0)}{\phi(-t)\theta(0,t)},
\end{equation}
where $i(-t,0)$ is the number of initially infected individuals.

The number of infected individuals who are infected at time $t$ is the sum of all individuals who have become infected at some prior time and not yet recovered. This can be split into the individuals who were initially infected and remain infected at time $t$, represented by $I_0(t)$, and the sum of individuals who have been infected at some $t>0$. Hence, we can write the number of infected individuals at time $t$ as
\begin{equation}\label{eq:full_I}
I(t)=\int_0^t \theta(t,t')\phi(t-t') q^+(I,t')dt'+I_0(t),
\end{equation}
where $I_0(t)$ can be represented as
\begin{equation}\label{eq:I_0}
I_0(t)=\int_{-\infty}^0 \frac{\theta(t,t')\phi(t-t')}{\theta(0,t')\phi(-t')} i(-t',0)dt'.
\end{equation}
We will take the initial conditions of the infected population to be
\begin{equation}\label{eq_delta} 
i(-t,0)=i_0\delta(-t),
\end{equation}
where $i_0$ is a constant. This simplifies Eq. \eqref{eq:I_0} to,
 \begin{equation} \label{eq_I0}
 I_0(t)=\theta(t,0)\phi(t) i_0.
 \end{equation}
Taking the derivative of Eq. \eqref{eq:full_I} with the initial conditions defined by Eqs. \eqref{eq_delta} and \eqref{eq_I0}, we arrive at
\begin{equation} \label{eq:dIdt}
\begin{split}
\frac{dI(t)}{dt}=&\omega(t)S(t)\left(\int_{0}^t \rho(t-t')\theta(t,t')\phi(t-t')q^+(I,t')dt'+\rho(t)\phi(t,0)\theta(t,0)i_0\right)\\&-\int_0^t \theta(t,t')\psi(t-t') q^+(I,t')dt' 
-\theta(t,0)\psi(t)i_0- \gamma(t)I(t).
\end{split}
\end{equation}
In order to write the derivative in terms of $I(t)$ instead of $q^+(I,t)$, we first define
\begin{equation}
F_I(t)=\int_{0}^t \rho(t-t')\theta(t,t')\phi(t-t')q^+(I,t')dt'
\end{equation}
and 
\begin{equation}
F_R(t)=\int_0^t \theta(t,t')\psi(t-t') q^+(I,t')dt. 
\end{equation}
Using the semi-group property and Laplace transforms, $F_I(t)$ and $F_R(t)$ can be written respectively as,
\begin{align}
\mathcal{L}_t\left\{\frac{F_I(t)}{\theta(t,0)}\right\}&=\mathcal{L}_t\left\{\rho(t)\phi(t)\right\}\mathcal{L}\left\{\frac{q^+(I,t)}{\theta(t,0)}\right\},\\
\mathcal{L}_t\left\{\frac{F_R(t)}{\theta(t,0)}\right\}&=\mathcal{L}_t\left\{\psi(t)\right\}\mathcal{L}_t\left\{\frac{q^+(I,t)}{\theta(t,0)}\right\}.
\end{align}
The Laplace transform of Eq. \eqref{eq:full_I} is
\begin{equation}\label{eq:lap_I}
\mathcal{L}_t\left\{\frac{I(t)-I_0(t)}{\theta(t,0)}\right\}=\mathcal{L}_t\left\{\phi(t)\right\}\mathcal{L}_t\left\{\frac{q^+(I,t)}{\theta(t,0)}\right\},
\end{equation}
from which we can write $F_I(t)$ as,
\begin{equation}\label{eq:F_inf}
F_I(t)=\theta(t,0)\int_0^t K_I(t-t')\frac{I(t')-I_0(t)}{\theta(t',0)}dt'
\end{equation}
where the infectivity kernel, $K_I$, is defined as
\begin{equation}\label{eq:inf_kernel}
K_I(t)=\mathcal{L}_s^{-1}\left\{\frac{\mathcal{L}_t\left\{\rho(t)\phi(t)\right\}}{\mathcal{L}_t\left\{\phi(t)\right\}}\right\}.
\end{equation}
Similarly, we can write $F_R(t)$ as
\begin{equation} \label{eq:F_rec}
F_R(t)=\theta(t,0)\int_0^t K_R(t-t')\frac{I(t')-I_0(t)}{\theta(t',0)}dt'
\end{equation}
where the recovery kernel, $K_R$, is
\begin{equation}\label{eq:rec_kernel}
K_R(t)=\mathcal{L}_s^{-1}\left\{\frac{\mathcal{L}_t\left\{\psi(t)\right\}}{\mathcal{L}_t\left\{\phi(t)\right\}}\right\}.
\end{equation}
Substituting Eqs. \eqref{eq:F_inf} and \eqref{eq:F_rec} into Eq. \eqref{eq:dIdt}, the governing equation for the infective compartment becomes
\begin{equation}
\frac{dI(t)}{dt}=\omega(t)S(t)\theta(t,0)\int_0^t K_I(t-t')\frac{I(t')}{\theta(t',0)}dt'-\theta(t,0)\int_0^t K_R(t-t')\frac{I(t')}{\theta(t',0)}dt'-\gamma(t)I(t).
\end{equation}
By considering the balance of flux between the compartments and the vital dynamics, we can write the governing equations for the Susceptible and Recovery compartments. Hence, the governing equation for the Susceptible compartment is,
\begin{equation}
\frac{dS(t)}{dt}=\lambda(t)- \omega(t)S(t)\theta(t,0)\int_0^t K_I(t-t')\frac{I(t')}{\theta(t',0)}dt'   -\gamma(t)S(t),
\end{equation}
and the governing equation for the Recovery compartment is,
\begin{equation}
\frac{dR(t)}{dt}=\theta(t,0)\int_0^t K_R(t-t')\frac{I(t')}{\theta(t',0)}dt'-\gamma(t)R(t).
\end{equation}
This gives us the full set of master equations for the time evolution of the epidemic across the Susceptible, Infective and Recovered populations.
%\section{Appendix} \label{sec_broader_initial_cond}
\section*{Acknowledgements}
This research was funded by Australian Research Council grant number DP200100345.

\bibliography{didrSIR.bib}


\begin{thebibliography}{29}
\ifx \bisbn   \undefined \def \bisbn  #1{ISBN #1}\fi
\ifx \binits  \undefined \def \binits#1{#1}\fi
\ifx \bauthor  \undefined \def \bauthor#1{#1}\fi
\ifx \batitle  \undefined \def \batitle#1{#1}\fi
\ifx \bjtitle  \undefined \def \bjtitle#1{#1}\fi
\ifx \bvolume  \undefined \def \bvolume#1{\textbf{#1}}\fi
\ifx \byear  \undefined \def \byear#1{#1}\fi
\ifx \bissue  \undefined \def \bissue#1{#1}\fi
\ifx \bfpage  \undefined \def \bfpage#1{#1}\fi
\ifx \blpage  \undefined \def \blpage #1{#1}\fi
\ifx \burl  \undefined \def \burl#1{\textsf{#1}}\fi
\ifx \doiurl  \undefined \def \doiurl#1{\url{https://doi.org/#1}}\fi
\ifx \betal  \undefined \def \betal{\textit{et al.}}\fi
\ifx \binstitute  \undefined \def \binstitute#1{#1}\fi
\ifx \binstitutionaled  \undefined \def \binstitutionaled#1{#1}\fi
\ifx \bctitle  \undefined \def \bctitle#1{#1}\fi
\ifx \beditor  \undefined \def \beditor#1{#1}\fi
\ifx \bpublisher  \undefined \def \bpublisher#1{#1}\fi
\ifx \bbtitle  \undefined \def \bbtitle#1{#1}\fi
\ifx \bedition  \undefined \def \bedition#1{#1}\fi
\ifx \bseriesno  \undefined \def \bseriesno#1{#1}\fi
\ifx \blocation  \undefined \def \blocation#1{#1}\fi
\ifx \bsertitle  \undefined \def \bsertitle#1{#1}\fi
\ifx \bsnm \undefined \def \bsnm#1{#1}\fi
\ifx \bsuffix \undefined \def \bsuffix#1{#1}\fi
\ifx \bparticle \undefined \def \bparticle#1{#1}\fi
\ifx \barticle \undefined \def \barticle#1{#1}\fi
\bibcommenthead
\ifx \bconfdate \undefined \def \bconfdate #1{#1}\fi
\ifx \botherref \undefined \def \botherref #1{#1}\fi
\ifx \url \undefined \def \url#1{\textsf{#1}}\fi
\ifx \bchapter \undefined \def \bchapter#1{#1}\fi
\ifx \bbook \undefined \def \bbook#1{#1}\fi
\ifx \bcomment \undefined \def \bcomment#1{#1}\fi
\ifx \oauthor \undefined \def \oauthor#1{#1}\fi
\ifx \citeauthoryear \undefined \def \citeauthoryear#1{#1}\fi
\ifx \endbibitem  \undefined \def \endbibitem {}\fi
\ifx \bconflocation  \undefined \def \bconflocation#1{#1}\fi
\ifx \arxivurl  \undefined \def \arxivurl#1{\textsf{#1}}\fi
\csname PreBibitemsHook\endcsname

%%% 1
\bibitem[\protect\citeauthoryear{Kermack and McKendrick}{1927}]{KM1927}
\begin{barticle}
\bauthor{\bsnm{Kermack}, \binits{W.O.}},
\bauthor{\bsnm{McKendrick}, \binits{A.G.}}:
\batitle{A contribution to the mathematical theory of epidemics}.
\bjtitle{Proceedings of the Royal Society of London Series A}
\bvolume{115}(\bissue{772}),
\bfpage{700}--\blpage{721}
(\byear{1927})
\doiurl{10.1098/rspa.1927.0118}
\end{barticle}
\endbibitem

%%% 2
\bibitem[\protect\citeauthoryear{Kermack and McKendrick}{1932}]{KM1932}
\begin{barticle}
\bauthor{\bsnm{Kermack}, \binits{W.O.}},
\bauthor{\bsnm{McKendrick}, \binits{A.G.}}:
\batitle{Contributions to the mathematical theory of epidemics. {II}. the
  problem of endemicity}.
\bjtitle{Proceedings of the Royal Society of London Series A}
\bvolume{138}(\bissue{834}),
\bfpage{55}--\blpage{83}
(\byear{1932})
\doiurl{10.1098/rspa.1932.0171}
\end{barticle}
\endbibitem

%%% 3
\bibitem[\protect\citeauthoryear{Hethcote}{2000}]{H2000}
\begin{barticle}
\bauthor{\bsnm{Hethcote}, \binits{H.W.}}:
\batitle{The mathematics of infectious diseases}.
\bjtitle{SIAM Review}
\bvolume{42}(\bissue{4}),
\bfpage{599}--\blpage{653}
(\byear{2000})
\doiurl{10.1137/S0036144500371907}
\end{barticle}
\endbibitem

%%% 4
\bibitem[\protect\citeauthoryear{Murray}{2003}]{M2003}
\begin{bbook}
\bauthor{\bsnm{Murray}, \binits{J.D.}}:
\bbtitle{Mathematical Biology: {II}: {S}patial Models and Biomedical
  Applications}
vol. \bseriesno{18}.
\bpublisher{Springer},
\blocation{New York}
(\byear{2003})
\end{bbook}
\endbibitem

%%% 5
\bibitem[\protect\citeauthoryear{Hethcote and Tudor}{1980}]{HT1980}
\begin{barticle}
\bauthor{\bsnm{Hethcote}, \binits{H.W.}},
\bauthor{\bsnm{Tudor}, \binits{D.W.}}:
\batitle{Integral equation models for endemic infectious diseases}.
\bjtitle{Journal of mathematical biology}
\bvolume{9},
\bfpage{37}--\blpage{47}
(\byear{1980})
\doiurl{10.1007/BF00276034}
\end{barticle}
\endbibitem

%%% 6
\bibitem[\protect\citeauthoryear{Robertson et~al.}{2018}]{RHRC2018}
\begin{barticle}
\bauthor{\bsnm{Robertson}, \binits{S.L.}},
\bauthor{\bsnm{Henson}, \binits{S.M.}},
\bauthor{\bsnm{Robertson}, \binits{T.}},
\bauthor{\bsnm{Cushing}, \binits{J.M.}}:
\batitle{A matter of maturity: To delay or not to delay? {C}ontinuous-time
  compartmental models of structured populations in the literature 2000--2016}.
\bjtitle{Natural Resource Modeling}
\bvolume{31}(\bissue{1}),
\bfpage{12160}
(\byear{2018})
\doiurl{10.1111/nrm.12160}
\end{barticle}
\endbibitem

%%% 7
\bibitem[\protect\citeauthoryear{Cooke}{1979}]{C1979}
\begin{barticle}
\bauthor{\bsnm{Cooke}, \binits{K.L.}}:
\batitle{Stability analysis for a vector disease model}.
\bjtitle{The Rocky Mountain Journal of Mathematics}
\bvolume{9}(\bissue{1}),
\bfpage{31}--\blpage{42}
(\byear{1979})
\doiurl{10.1216/RMJ-1979-9-1-31}
\end{barticle}
\endbibitem

%%% 8
\bibitem[\protect\citeauthoryear{Cooke and {Van den Driessche}}{1996}]{CD1996}
\begin{barticle}
\bauthor{\bsnm{Cooke}, \binits{K.L.}},
\bauthor{\bsnm{{Van den Driessche}}, \binits{P.}}:
\batitle{Analysis of an {SEIRS} epidemic model with two delays}.
\bjtitle{Journal of mathematical biology}
\bvolume{35},
\bfpage{240}--\blpage{260}
(\byear{1996})
\doiurl{10.1007/s002850050051}
\end{barticle}
\endbibitem

%%% 9
\bibitem[\protect\citeauthoryear{Tchuenche et~al.}{2007}]{TNL2007}
\begin{barticle}
\bauthor{\bsnm{Tchuenche}, \binits{J.M.}},
\bauthor{\bsnm{Nwagwo}, \binits{A.}},
\bauthor{\bsnm{Levins}, \binits{R.}}:
\batitle{Global behaviour of an {SIR} epidemic model with time delay}.
\bjtitle{Mathematical Methods in the Applied Sciences}
\bvolume{30}(\bissue{6}),
\bfpage{733}--\blpage{749}
(\byear{2007})
\doiurl{10.1002/mma.810}
\end{barticle}
\endbibitem

%%% 10
\bibitem[\protect\citeauthoryear{Ma et~al.}{2004}]{MST2004}
\begin{barticle}
\bauthor{\bsnm{Ma}, \binits{W.}},
\bauthor{\bsnm{Song}, \binits{M.}},
\bauthor{\bsnm{Takeuchi}, \binits{Y.}}:
\batitle{Global stability of an {SIR} epidemic model with time delay}.
\bjtitle{Applied Mathematics Letters}
\bvolume{17}(\bissue{10}),
\bfpage{1141}--\blpage{1145}
(\byear{2004})
\doiurl{10.1016/j.aml.2003.11.005}
\end{barticle}
\endbibitem

%%% 11
\bibitem[\protect\citeauthoryear{Beretta and Takeuchi}{1995}]{BT1995}
\begin{barticle}
\bauthor{\bsnm{Beretta}, \binits{E.}},
\bauthor{\bsnm{Takeuchi}, \binits{Y.}}:
\batitle{Global stability of an {SIR} epidemic model with time delays}.
\bjtitle{Journal of Mathematical Biology}
\bvolume{33}(\bissue{3}),
\bfpage{250}--\blpage{260}
(\byear{1995})
\doiurl{10.1007/BF00169563}
\end{barticle}
\endbibitem

%%% 12
\bibitem[\protect\citeauthoryear{Hethcote and {Van den
  Driessche}}{2000}]{HD2000}
\begin{barticle}
\bauthor{\bsnm{Hethcote}, \binits{H.W.}},
\bauthor{\bsnm{{Van den Driessche}}, \binits{P.}}:
\batitle{Two {SIS} epidemiologic models with delays}.
\bjtitle{Journal of Mathematical Biology}
\bvolume{40},
\bfpage{3}--\blpage{26}
(\byear{2000})
\doiurl{10.1007/s002850050003}
\end{barticle}
\endbibitem

%%% 13
\bibitem[\protect\citeauthoryear{Liu}{2015}]{Liu2015}
\begin{barticle}
\bauthor{\bsnm{Liu}, \binits{L.}}:
\batitle{A delayed {SIR} model with general nonlinear incidence rate}.
\bjtitle{Advances in Difference Equations}
\bvolume{2015}(\bissue{1}),
\bfpage{329}
(\byear{2015})
\doiurl{10.1186/s13662-015-0619-z}
\end{barticle}
\endbibitem

%%% 14
\bibitem[\protect\citeauthoryear{Yan and Liu}{2006}]{YL2006}
\begin{barticle}
\bauthor{\bsnm{Yan}, \binits{P.}},
\bauthor{\bsnm{Liu}, \binits{S.}}:
\batitle{{SEIR} epidemic model with delay}.
\bjtitle{The ANZIAM Journal}
\bvolume{48}(\bissue{1}),
\bfpage{119}--\blpage{134}
(\byear{2006})
\doiurl{10.1017/S144618110000345X}
\end{barticle}
\endbibitem

%%% 15
\bibitem[\protect\citeauthoryear{Huang et~al.}{2010}]{HTMW2010}
\begin{barticle}
\bauthor{\bsnm{Huang}, \binits{G.}},
\bauthor{\bsnm{Takeuchi}, \binits{Y.}},
\bauthor{\bsnm{Ma}, \binits{W.}},
\bauthor{\bsnm{Wei}, \binits{D.}}:
\batitle{Global stability for delay {SIR} and {SEIR} epidemic models with
  nonlinear incidence rate}.
\bjtitle{Bulletin of Mathematical Biology}
\bvolume{72},
\bfpage{1192}--\blpage{1207}
(\byear{2010})
\doiurl{10.1007/s11538-009-9487-6}
\end{barticle}
\endbibitem

%%% 16
\bibitem[\protect\citeauthoryear{Montroll and Weiss}{1965}]{MW1965}
\begin{barticle}
\bauthor{\bsnm{Montroll}, \binits{E.W.}},
\bauthor{\bsnm{Weiss}, \binits{G.H.}}:
\batitle{Random walks on lattices. {II}}.
\bjtitle{Journal of Mathematical Physics}
\bvolume{6}(\bissue{2}),
\bfpage{167}--\blpage{181}
(\byear{1965})
\doiurl{10.1063/1.1704269}
\end{barticle}
\endbibitem

%%% 17
\bibitem[\protect\citeauthoryear{Angstmann et~al.}{2017}]{AHM2017}
\begin{barticle}
\bauthor{\bsnm{Angstmann}, \binits{C.N.}},
\bauthor{\bsnm{Henry}, \binits{B.I.}},
\bauthor{\bsnm{McGann}, \binits{A.V.}}:
\batitle{A fractional-order infectivity and recovery {SIR} model}.
\bjtitle{Fractal and Fractional}
\bvolume{1}(\bissue{1}),
\bfpage{11}
(\byear{2017})
\doiurl{10.3390/fractalfract1010011}
\end{barticle}
\endbibitem

%%% 18
\bibitem[\protect\citeauthoryear{Angstmann et~al.}{2016a}]{AHM2016}
\begin{barticle}
\bauthor{\bsnm{Angstmann}, \binits{C.N.}},
\bauthor{\bsnm{Henry}, \binits{B.I.}},
\bauthor{\bsnm{McGann}, \binits{A.V.}}:
\batitle{A fractional order recovery {SIR} model from a stochastic process}.
\bjtitle{Bulletin of Mathematical Biology}
\bvolume{78},
\bfpage{468}--\blpage{499}
(\byear{2016})
\doiurl{10.1007/s11538-016-0151-7}
\end{barticle}
\endbibitem

%%% 19
\bibitem[\protect\citeauthoryear{Angstmann et~al.}{2016b}]{AHM2016fi}
\begin{barticle}
\bauthor{\bsnm{Angstmann}, \binits{C.N.}},
\bauthor{\bsnm{Henry}, \binits{B.I.}},
\bauthor{\bsnm{McGann}, \binits{A.V.}}:
\batitle{A fractional-order infectivity {SIR} model}.
\bjtitle{Physica A: Statistical Mechanics and its Applications}
\bvolume{452},
\bfpage{86}--\blpage{93}
(\byear{2016})
\doiurl{10.1016/j.physa.2016.02.029}
\end{barticle}
\endbibitem

%%% 20
\bibitem[\protect\citeauthoryear{Angstmann et~al.}{2017}]{AEHMMN2017}
\begin{barticle}
\bauthor{\bsnm{Angstmann}, \binits{C.N.}},
\bauthor{\bsnm{Erickson}, \binits{A.M.}},
\bauthor{\bsnm{Henry}, \binits{B.I.}},
\bauthor{\bsnm{McGann}, \binits{A.V.}},
\bauthor{\bsnm{Murray}, \binits{J.M.}},
\bauthor{\bsnm{Nichols}, \binits{J.A.}}:
\batitle{Fractional order compartment models}.
\bjtitle{SIAM Journal on Applied Mathematics}
\bvolume{77}(\bissue{2}),
\bfpage{430}--\blpage{446}
(\byear{2017})
\doiurl{10.1137/16M1069249}
\end{barticle}
\endbibitem

%%% 21
\bibitem[\protect\citeauthoryear{Angstmann et~al.}{2021}]{AEHMMN2021}
\begin{barticle}
\bauthor{\bsnm{Angstmann}, \binits{C.N.}},
\bauthor{\bsnm{Erickson}, \binits{A.M.}},
\bauthor{\bsnm{Henry}, \binits{B.I.}},
\bauthor{\bsnm{McGann}, \binits{A.V.}},
\bauthor{\bsnm{Murray}, \binits{J.M.}},
\bauthor{\bsnm{Nichols}, \binits{J.A.}}:
\batitle{A general framework for fractional order compartment models}.
\bjtitle{SIAM Review}
\bvolume{63}(\bissue{2}),
\bfpage{375}--\blpage{392}
(\byear{2021})
\doiurl{10.1137/21M1398549}
\end{barticle}
\endbibitem

%%% 22
\bibitem[\protect\citeauthoryear{Angstmann et~al.}{2023}]{ABHJX2023}
\begin{barticle}
\bauthor{\bsnm{Angstmann}, \binits{C.N.}},
\bauthor{\bsnm{Burney}, \binits{S.-J.M.}},
\bauthor{\bsnm{Henry}, \binits{B.I.}},
\bauthor{\bsnm{Jacobs}, \binits{B.A.}},
\bauthor{\bsnm{Xu}, \binits{Z.}}:
\batitle{A systematic approach to delay functions}.
\bjtitle{Mathematics}
\bvolume{11}(\bissue{21}),
\bfpage{4526}
(\byear{2023})
\doiurl{10.3390/math11214526}
\end{barticle}
\endbibitem

%%% 23
\bibitem[\protect\citeauthoryear{Angstmann et~al.}{2024}]{ABHHMX2024}
\begin{barticle}
\bauthor{\bsnm{Angstmann}, \binits{C.N.}},
\bauthor{\bsnm{McGann}, \binits{A.V.}},
\bauthor{\bsnm{Xu}, \binits{Z.}}:
\batitle{Delay compartment models from a stochastic process}.
\bjtitle{arXiv}
\bvolume{[math.DS]}(\bissue{2406.17242}),
\bfpage{1}--\blpage{20}
(\byear{2024})
\doiurl{10.48550/arXiv.2406.17242}
\end{barticle}
\endbibitem

%%% 24
\bibitem[\protect\citeauthoryear{Couch and Kasel}{1983}]{CK1983}
\begin{barticle}
\bauthor{\bsnm{Couch}, \binits{R.B.}},
\bauthor{\bsnm{Kasel}, \binits{J.A.}}:
\batitle{Immunity to influenza in man.}
\bjtitle{Annual review of microbiology}
\bvolume{37},
\bfpage{529}--\blpage{549}
(\byear{1983})
\doiurl{10.1146/annurev.mi.37.100183.002525}
\end{barticle}
\endbibitem

%%% 25
\bibitem[\protect\citeauthoryear{Peteranderl et~al.}{2016}]{PHS2016}
\begin{bchapter}
\bauthor{\bsnm{Peteranderl}, \binits{C.}},
\bauthor{\bsnm{Herold}, \binits{S.}},
\bauthor{\bsnm{Schmoldt}, \binits{C.}}:
\bctitle{Human influenza virus infections}.
In: \bbtitle{Seminars in Respiratory and Critical Care Medicine},
vol. \bseriesno{37},
pp. \bfpage{487}--\blpage{500}
(\byear{2016}).
\doiurl{10.1055/s-0036-1584801} .
\bcomment{Thieme Medical Publishers}
\end{bchapter}
\endbibitem

%%% 26
\bibitem[\protect\citeauthoryear{Insinga et~al.}{2007}]{IDELB2007}
\begin{barticle}
\bauthor{\bsnm{Insinga}, \binits{R.P.}},
\bauthor{\bsnm{Dasbach}, \binits{E.J.}},
\bauthor{\bsnm{Elbasha}, \binits{E.H.}},
\bauthor{\bsnm{Liaw}, \binits{K.}},
\bauthor{\bsnm{Barr}, \binits{E.}}:
\batitle{Incidence and duration of cervical human papillomavirus 6, 11, 16, and
  18 infections in young women: an evaluation from multiple analytic
  perspectives}.
\bjtitle{Cancer Epidemiology Biomarkers \& Prevention}
\bvolume{16}(\bissue{4}),
\bfpage{709}--\blpage{715}
(\byear{2007})
\doiurl{10.1158/1055-9965.EPI-06-0846}
\end{barticle}
\endbibitem

%%% 27
\bibitem[\protect\citeauthoryear{Manfredi and D'Onofrio}{2013}]{MD2013}
\begin{bbook}
\bauthor{\bsnm{Manfredi}, \binits{P.}},
\bauthor{\bsnm{D'Onofrio}, \binits{A.}}:
\bbtitle{Modeling the Interplay Between Human Behavior and the Spread of
  Infectious Diseases}.
\bpublisher{Springer},
\blocation{New York}
(\byear{2013})
\end{bbook}
\endbibitem

%%% 28
\bibitem[\protect\citeauthoryear{Guiro et~al.}{2023}]{GOO2023}
\begin{bchapter}
\bauthor{\bsnm{Guiro}, \binits{A.}},
\bauthor{\bsnm{Ouedraogo}, \binits{D.}},
\bauthor{\bsnm{Ouedraogo}, \binits{H.}}:
\bctitle{Global stability for a delay {SIR} epidemic model with general
  incidence function, observers design}.
In: \bbtitle{Partial Differential Equations and Applications},
pp. \bfpage{259}--\blpage{280}
(\byear{2023}).
\doiurl{10.1007/978-3-031-27661-3_10} .
\bcomment{Springer Proceedings in Mathematics \& Statistics}
\end{bchapter}
\endbibitem

%%% 29
\bibitem[\protect\citeauthoryear{Gu}{2010}]{G2010}
\begin{barticle}
\bauthor{\bsnm{Gu}, \binits{K.}}:
\batitle{Stability problem of systems with multiple delay channels}.
\bjtitle{Automatica}
\bvolume{46}(\bissue{4}),
\bfpage{743}--\blpage{751}
(\byear{2010})
\doiurl{10.1016/j.automatica.2010.01.028}
\end{barticle}
\endbibitem

\end{thebibliography}
\end{document}